\documentstyle[12pt]{report}

\renewcommand\c{\circ}
\newcommand\C{{\mbox{\rm\bf C\hspace{-7.3pt}}{^{_{\bf\mid}}}\hspace{4.5pt}}}
\newcommand\Ca{\mbox{\boldmath $C\! a$}}

\newcommand\D{{\mbox{\rm\bf D\hspace{-7.3pt}}{^{_{\bf\mid}}}\hspace{4.5pt}}}

\newcommand\G{{\cal G}}

\newcommand\hsb{\hskip-14pt}
\newcommand\hsp{\hskip-24pt}
\renewcommand\H{{\rm I\hspace{-2.5pt} H}}

\newcommand\I{{\cal I}}

\newcommand\ok{{\cal O}}

\newcommand\op[1]{\mathop{\rm #1}\nolimits}

\newcommand\po{$\!\!\!{\bf .}$ }

\newcommand\R{{\rm I\hspace{-2.5pt} R}}

\newcommand\ve{\varepsilon}
\newcommand\vp{\varphi}

\def\Rom#1{\uppercase\expandafter{\romannumeral#1}}

\newcommand\1{{\bf 1}}

\newcommand\qed{\phantom{\underline{y}}\hfill\hfill$\Box$}

\newcommand\bib[1]{\bibitem[#1]{#1}}

\newcommand{\text}[1]{{\mbox{\rm #1}}}
\newcommand{\dfrac}[2]{\frac{\displaystyle #1}{\displaystyle #2}}

\newtheorem{th}{Theorem}
\newtheorem{prop}{Proposition}

{\endtrivlist}
\newenvironment{lem}{\trivlist \item[\hskip \labelsep{\bf Lemma.}]\it}%
{\endtrivlist}
\newenvironment{cor}{\trivlist \item[\hskip \labelsep{\bf Corollary.}]\it}%
{\endtrivlist}
\newenvironment{dfn}{\trivlist \item[\hskip \labelsep{\bf Definition.}]}%
{\endtrivlist}
{\endtrivlist}
\newenvironment{proof}{\trivlist \item[\hskip
\labelsep{{\it\underline{Proof}.\/}}]}%
{\endtrivlist}
{\endtrivlist}

\newcounter{a}
\newcounter{f}
\makeatletter
\newcommand{\@thefnmark}{$^\fnsymbol{f}$}
\renewcommand{\@makefnmark}{\hbox{\mathsurround=0pt
                           $^{\fnsymbol{f}}$}}
\renewcommand{\@makefntext}[1]{\parindent=1em\noindent
            \hbox to 1.8em{\hss$^{\fnsymbol{f}}$}#1}
\makeatother

\begin{document}

\title{The Kobayashi pseudodistance \\ on almost complex manifolds}

\author{Boris~S.~Kruglikov and Marius Overholt}

\date{March 1, 1997}

\begin{abstract}

We extend the definition of the Kobayashi pseudodistance to
almost complex manifolds and show that its familiar properties are for the
most part preserved.  We also study the automorphism group
of an almost complex manifold and finish with some examples.

\end{abstract}

\maketitle

\tableofcontents

\clearpage


\chapter*{Introduction}
\addcontentsline{toc}{chapter}{\bf\quad \  Introduction}

\hspace{13.5pt}
The Poincar\'{e} metric on the open unit disk $\D$
in the complex plane $\C$ is a Riemannian metric
 $$
|v|=\frac{|v|_{euc}}{1-|z|^2}
 $$
conformal with the Euclidean metric $|\cdot|_{euc}$, that induces a
distance $d$ on $\D$ with the remarkable
property that every holomorphic mapping $f : \D \to \D$ is distance
nonincreasing in $d$. This fact, discovered in 1915 by Pick \cite{Pi} is
an invariant formulation of the Schwarz lemma.  In 1967 Kobayashi
used the distance $d$ to define a pseudodistance $d_M$ on any (connected)
complex manifold $M$, such that any
holomorphic mapping from a complex manifold $L$ to a complex manifold $M$
is distance nonincreasing with respect to $d_L$ and $d_M$.  When this
Kobayashi pseudodistance $d_M$ is a distance, it may be used to obtain
information about holomorphic mappings to or from $M$; in this situation
$M$ is said to be (Kobayashi-) hyperbolic.  Some references for Kobayashi
hyperbolicity are \cite{JP,Ko1,Ko2,La,NO}.

In this paper\footnote{This paper appeared during the visit of
the first
author to the University of Tromsoe. He thanks professor
V.\,Lychagin who had invited and greeted him and other
mathematicians from the Department for hospitality and kind
attention.} we extend
the definition of the Kobayashi pseudodistance to
almost complex manifolds and show that its familiar properties are for the
most part preserved.  Results whose proofs are similar to known ones for
complex manifolds are merely stated.  We also study the automorphism group
of an almost complex manifold and finish with some examples.

\chapter{Definition of the pseudometric}
\markboth{}{}

\hspace{13.5pt}
Let  $\C$ denote the complex plane with its standard complex structure,
$\D$ the open unit disk in it and $e=1\in T_0\D$.
Let $(M^{2n},J)$ be an almost complex manifold, which in this paper
will be taken to be connected and $C^\infty$.  This means that the smooth
field of automorphisms  $J\in T^*M \otimes TM$  satisfies  $J^2 = -\1$.
A Kobayashi chain joining two
points $p,q\in M$ is a sequence of pseudoholomorphic mappings
 $$
f_k : \D \rightarrow (M^{2n},J), \ \ \ k=1,...,m,
 $$
and points $z_k,w_k \in \D$ such that $f_1(z_1) = p, \ f_m(w_m) = q$
and $f_{k}(w_k) = f_{k+1}(z_{k+1})$.  The {\it Kobayashi pseudodistance\/}
from $p$ to $q$ is defined by
 $$
d_M(p,q) = \op{inf} \ \sum_{k=1}^{m}d(z_k,w_k)
 $$
where the infimum is taken over all Kobayashi chains joining $p$ to $q$,
if there exists some Kobayashi chain joining $p$ to $q$, and is defined to
be $+\infty$ otherwise.  It is well known that $d_{\C} \equiv 0$ and $d_{\D}
 \equiv d$.

 \begin{prop}\po
The function $d_M : M \times M \to \R$ is nonnegative, symmetric
and satisfies the triangle inequality, i.e.\ it is a pseudodistance, called
the Kobayashi pseudodistance.  Pseudoholomorphic mappings
between almost complex manifolds are distance nonincreasing with respect
to the Kobayashi pseudodistance.
 \end{prop}

As in the case of complex manifolds, $d_M$ is finite on any almost complex
manifold, but this requires an existence theorem for pseudoholomorphic
disks due to Nijenhuis and Woolf \cite{NW}.

 \begin{th}\po
{\rm(i)}. The Kobayashi pseudodistance $d_M$ is finite and continuous on $M
\times M$
for any almost complex manifold $(M^{2n},J)$.\newline
{\rm(ii)}. If $d_M$ is a distance, it induces the standard topology on $M$.
 \end{th}

 \begin{proof}
We must first prove that any two points in $M$ can be joined by a
Kobayashi chain.  It is enough to prove that for any point
$p \in M$ there is a neighborhood $U_p$ of $p$ such that
any point $q$ in $U_p$ can be joined to $p$ by a single
pseudoholomorphic disk.  For then the set of points that can
be joined by a Kobayashi chain to some fixed point is open
and nonempty, and its complement is open,
so every point can be joined to the fixed point since $M$
is connected.

The problem is now local, so we consider $({\R}^{2n},J)$.
Let $v \in T_{0}{\R}^{2n} = {\R}^{2n}$.
By theorem III of \cite{NW}, there is
some neighborhood $V$ of $0$ in ${\R}^{2n}$ such that
if $v \in V$, then there exists a pseudoholomorphic mapping
$f : {\D} \rightarrow ({\R}^{2n},J)$ with $f(0) = 0$ and
$f_{\ast}(0)e = v$; this mapping could be choosen canonically.
We write this mapping as $f(z;v)$ to show its
dependence on $v$. We denote its differential with respect to $z$
by $f_{\ast}(z;v)$. By 5.4a
of \cite{NW}, $f : {\D} \times V \rightarrow
{\R}^{2n}$ is a $C^\infty$ mapping.  It satisfies
$f(0;v) = 0$ and $f_{*}(0;v)e = v$.  In addition, by 5.2a of \cite{NW}
, it satisfies $f(z;0) = 0$ on the unit disk.  Since
 $$
\frac{\partial f_{*}(0;v)e}{\partial v} = I
 $$
there is some neighborhood $W$ of $0$ in $\D$ such that for every
$\zeta \in W$, the matrix
 $$
\left.\frac{\partial f_{*}(\zeta;v)e}{\partial v}\right|_{v=0}
 $$
is of full rank.  Now
 $$
\left.\frac{\partial f(\zeta;v)}{\partial v}\right|_{v=0} =
\left.{\zeta}\frac{\partial f_{*}(\zeta;v)e}{\partial v}\right|_{v=0} +
O(|\zeta|^2),
 $$
where multiplication is with respect to the complex structure on ${\R}^{2n}
\simeq \C^n$ induced by the linearization.
So there exists some $\zeta_0 \in W \backslash {\{}0{\}}$ such that
 $$
\left.\frac{\partial f(\zeta_0;v)}{\partial v}\right|_{v=0}
 $$
is a matrix of full rank.  Then $f(\zeta_0;V)$ contains an open
neighborhood of $0$ in ${\R}^{2n}$, because $f(\zeta_0;0) = 0$.  So if
$q$ is in this neighborhood, there is some $v \in V$ such that
$f(\cdot;v)$ is a pseudoholomorphic disk joining $0$ and $q$.

Since
 $$
|d_M(p,q) - d_M(p_0,q_0)| \leq d_M(p,p_0) + d_M(q,q_0)
 $$
by the triangle inequality, to prove continuity of $d_M$ it is enough to
prove
continuity of the mappings $d(\cdot,p_0)$ for all $p_0 \in M$.  The
reasoning
above shows that
for every $\ve > 0$ there exists some neighborhood $U$ of $p_0$ such that
for
some $\zeta_0 \in \D$ with $d(0,\zeta_0) = \ve$ the image $f(\zeta_0;V)$
contains  $U$  and thus  $d_M(p,p_0) \leq \ve$  for any  $p \in U$.
The proof that $d_M$ induces the standard
topology on $M$ if it is a distance can now be carried through by the
method
of Barth \cite{Ba}.\qed
 \end{proof}

It is obvious that for any  $p_1,q_1 \in (M_1,J_1)$  and
 $p_2,q_2 \in (M_2,J_2)$
we have
 $$
  \begin{array}{l}
\max\, [d_{M_1}(p_1,q_1),d_{M_2}(p_2,q_2)] \leq \\
\qquad\qquad\qquad\quad\
d_{M_1 \times M_2}(p_1\times p_2,q_1\times q_2) \leq
d_{M_1}(p_1,q_1) + d_{M_2}(p_2,q_2),
  \end{array}
 $$
but this inequality is much less useful than in the complex case.  For if
$({M_1} \times {M_2},J)$ is a product manifold with
an almost complex structure $J$, it is almost never true even locally that
$J$
decomposes as a product of complex structures $J = J_1 \times J_2$.


In chapter 2 we state a result on hyperbolicity of almost complex manifolds
which are fiber spaces.  A more elementary fact is

 \begin{prop}\po
Let  $\pi : E \rightarrow B$ be an almost complex locally trivial bundle.
Then for every two points $p,q \in B$ we have
 $$
d_{B}(p,q) = d_{E}({\pi}^{-1}(p),{\pi}^{-1}(q)) =
\inf {\{}d_E(\tilde p,\tilde q) \,|\, \pi(\tilde{p})=p ,
\pi(\tilde{q}) =q{\}}.
 $$
 \end{prop}

 \begin{proof}
Since $\pi : E \to B$ is pseudoholomorphic we have
 $$
d_B(p,q) \leq d_E(\tilde{p},\tilde{q}).
 $$
On the other hand every pseudoholomorphic disk $f : \D \to B$
lifts to a pseudoholomorphic disk $\tilde{f} : \D \to E$.\qed
 \end{proof}

 \begin{cor}
The same formula holds for almost complex (unbranched) coverings.
Hyperbolicity of an almost complex manifold is equivalent to
hyperbolicity of any of its (unbranched) coverings.
 \end{cor}

We note that the Kobayashi-Royden infinitesimal pseudometric  $F_M$  may
also
be defined on almost complex manifolds using the results of Nijenhuis and
Woolf.  The definition is
 $$
F_M(p;v) = \inf \frac{1}{r}
 $$
where  $p \in M$, $v \in T_{p}M$  and the infimum is taken with respect to
all
pseudoholomorphic mappings  $f : \D \rightarrow M$ with $f(0) = p$  and
$f_{\ast}(0)e = rv$.  This may be reformulated as another invariant
${\pi}F_M(p;v)^{-2} =\sup({\pi}r^2)$ which looks similar to the so-called
Gromov
width from symplectic geometry, see \cite{Gr,MS}.

\chapter{Hyperbolicity}
\markboth{}{}

\hspace{13.5pt}
An almost complex manifold $(M^{2n},J)$ is {\it Kobayashi hyperbolic\/}
if the Kobayashi pseudodistance $d_M$ is a distance.  It is
{\it Brody hyperbolic\/} if any pseudoholomorphic mapping
$f : {\C} \to M$ is constant. It is clear that Kobayashi hyperbolicity
implies Brody hyperbolicity. Brody \cite{Br} discovered that for complex
manifolds the
converse holds under the additional condition of compactness. As is
usual, we take the term hyperbolic to mean Kobayashi hyperbolic,
and denote other hyperbolic properties by a prefix.

The following is a version of Brody's theorem \cite{Br} for almost complex
manifolds:

 \begin{th}\po
Let $(M^{2n},J)$ be an almost complex manifold and $|\cdot|$ a
continuous norm on $TM$.  If there exists a constant $C$ such that
 $$
|f_{\ast}(0)e| \leq C
 $$
for all pseudoholomorphic mappings $f : \D \to M$, then $M$ is
hyperbolic.

If $M$ is compact hyperbolic, then there exists a
constant $C$ such that
 $$
|f_{\ast}(z)e| \leq \frac{C}{1 - |z|^2}
 $$
for all pseudoholomorphic mappings $f : \D \to M$ and all $z \in \D$.

A compact almost complex manifold is hyperbolic if and only if it
is Brody hyperbolic.
 \end{th}

 \begin{proof}
Suppose that a constant $C$ exists with $|f_{\ast}(0)e| \leq C$ for
all pseudoholomorphic mappings $f : \D \to M$. Then
 $$
|f_{\ast}(z)e| \leq \frac{C}{1 - |z|^2}
 $$
by precomposing with suitable disk automorphisms.
Let $\rho$ be the distance associated to $|\cdot|$. Then for every
pseudoholomorphic mapping $f:\D\to M$ and for every pair $z_1, z_2\in
\D$ the following holds: $\rho(f(z_1),f(z_2))\le C d(z_1,z_2)$. The last
inequality implies $0< \rho(p,q)\le C d_M(p,q)$ for every pair of points
$p\ne q$, whence the first claim follows.

To prove the second claim, note that if no such constant $C$ exists,
then neither does there exist a constant $C$ such that
 $$
|f_{\ast}(0)e| \leq C
 $$
for all pseudoholomorphic mappings $f : \D \to M$, again by precomposing
with suitable disk automorphisms. So let $f_k : \D \to M$ be a
sequence of pseudoholomorphic mappings with
 $$
|(f_k)_{\ast}(0)e| \to \infty,
 $$
and assume that $M$ is compact.  We will show that $M$ is not hyperbolic.
Since $M$ is compact, we may extract a subsequence to assure that
$f_k(0) \to p \in M$.  Let $K$ be a compact set with $p$ an interior point
on which
$J$ is tamed by an exact symplectic form, and let $r_k$ be the supremum
of radii $r \leq 1$ such that $f_k(r\D) \subseteq K$.
Gromov's Schwarz lemma, see Corollary 4.1.4 of \cite{Mu}, implies that
$r_k \to 0$,
since $|(f_k)_{\ast}(0)e| \to \infty$.  It is clear that there exists
a sequence $z_k \in r_k\overline{\D}$ such that
$q_k = f(z_k) \in {\partial}K$.
Extract a subsequence to assure that $q_k \to q \in {\partial}K$.
Then $d_M(p,q_k) \to d_M(p,q)$ by the continuity of $d_M$ and
 $$
d_M(p,q_k) = d_M(f_k(0),f_k(z_k)) \leq d(0,z_k) \to 0
 $$
and so $d_M(p,q) = 0$, hence $M$ is not hyperbolic.

It is clear that a hyperbolic almost complex manifold is Brody hyperbolic.
Assume on the other hand that $M$ is a compact, but not hyperbolic,
almost complex manifold.  We shall prove that there exists a
nonconstant pseudoholomorphic mapping $g : \C \to M$.

We put a Riemannian metric $|\cdot|$ on $M$.  Since $M$ is not
hyperbolic, there is a sequence of pseudoholomorphic mappings
$f_k : \D \to M$ with $|(f_k)_{\ast}(0)e| \to \infty$.  The proof of
Brody's reparametrization lemma \cite{Br} goes through unchanged for almost
complex manifolds.  We use it in the form that if $(M^{2n},J)$ is an almost
complex manifold, $|\cdot|$ a continuous norm on $TM$ and $f : r\D \to M$
a pseudoholomorphic mapping such that $|f_{\ast}(0)e| \geq c > 0$, then
there exists a pseudoholomorphic mapping $h : r\D \to M$ such that
$|h_{\ast}(0)e| = c/2$,
 $$
|h_{\ast}(z)e| \leq \frac{c}{2}\frac{r^2}{r^2 - |z|^2},
 $$
and $h(r\D) \subseteq f(r\D)$.  See p. 27 of \cite{NO} for a proof that
can be adapted virtually unchanged.  We apply Brody's reparametrization
lemma to the sequence of pseudoholomorphic mappings
 $$
f_k\left(\frac{2z}{|(f_k)_{\ast}(0)e|}\right)
 $$
from $r_k\D$ to $M$, where $r_k = |(f_k)_{\ast}(0)e|/2 \to \infty$.
This gives a sequence of pseudoholomorphic mappings $g_k : r_k\D \to M$
with $|(g_k)_{\ast}(0)e| = 1$ and
 $$
|(g_k)_{\ast}(z)e| \leq \frac{r_k^2}{r_k^2 - |z|^2}.
 $$
Using the last inequality and compactness theorem 4.1.3 of \cite{MS},
we conclude that
$g_k$ has a subsequence which converges uniformly with all derivatives
on compact subsets of $\C$ to a pseudoholomorphic mapping $g : \C \to M$.
It is nonconstant because $|(g)_{\ast}(0)e| = 1$.\qed
\end{proof}

More general versions of Brody's theorem are known for complex manifolds.
See theorems 2.2 and 2.3 of \cite{La}.  They have analogs for almost
complex manifolds, which are easily proved by modifying the proofs in
\cite{La} along the lines of the above proof.
Then one can prove the following theorem in the same way as theorem 3.1
is proved in \cite{La}:

\begin{th}\po
Let $\pi : E \to B$ be a proper pseudoholomorphic mapping between
almost complex manifolds, whose fibers are manifolds.  If $B$ is
hyperbolic and each fiber is hyperbolic, then $E$ is hyperbolic.
If the fiber above a point of $B$ is hyperbolic, this point has a
neighborhood such that the fibers above all points of this neighborhood
are hyperbolic.
 \end{th}

\chapter{The automorphism group}
\markboth{}{}

\hspace{13.5pt}
In contrast to the complex case, generic almost complex
structures admit no local automorphisms other than the identity
and so the group $\op{Aut}(M^{2n},J)$ is usually trivial.
The obstruction theory for local and formal automorphisms of almost
complex manifolds was developed in \cite{Kr}. However in practice
we usually need criteria for this generic property. Here we give
one connected with the notion of hyperbolicity.

It is a result of Kobayashi \cite{Ko1} that the automorphism group
of a closed (usually called in the literature just compact but with
assumption of empty boundary) hyperbolic complex manifold is finite.
This is also true in the almost complex case.
For closed almost complex manifolds, the group
of pseudoholomorphic diffeomorphisms of the manifold to itself is a Lie
transformation group, when equipped with the topology of
$\Sigma$-convergence, i.e.\ uniform convergence on compact sets
of the mappings and their derivatives through the third order,
as was proved by Boothby, Kobayashi and Wang \cite{BKW}.


  \begin{th}\po
The automorphism group  $G = \op{Aut}(M^{2n},J)$  of a closed hyperbolic
almost complex manifold $(M^{2n},J)$ is finite.
 \end{th}

 \begin{proof}
Since  $G$  is a Lie group we may consider its Lie algebra $\G$ consisting
of pseudoholomorphic vector fields; i.e.\ fields $\xi$ on $M^{2n}$
such that  $L_{\xi}J = 0$ or  $[\xi,J\eta] = J[\xi,\eta]$ for every vector
field $\eta$ on $M^{2n}$.  In particular, $[\xi,J\xi] = 0$.  This means
that if $0 \neq \xi \in \G$ then we have a nonconstant pseudoholomorphic
curve $f : \C \rightarrow M^{2n}$ through a point $p \in M^{2n}$ of the
form
 $$
z = x + iy \mapsto \exp(x\xi + y(J\xi))p.
 $$
The exponential mapping is globally defined since $M^{2n}$ was assumed
closed.  Now $(M^{2n},J)$ was assumed hyperbolic and so such a
nonconstant pseudoholomorphic curve cannot exist.  Hence $\G$ is trivial
and thus $G$ is discrete.

We must now show that $G$ is compact.  Since $M$ is compact and second
countable, $G$ is second countable.  So it is enough to prove that any
sequence $\vp_k \in G$ has a $\Sigma$-convergent subsequence.  In the
following, each time we extract a subsequence, we give it the same
notation as the original sequence.  Note that
$G$ is a subgroup of the isometry group $\I(M,d_M)$ of $M$,
where $d_M$ is the Kobayashi distance on $M$.
According to the theorem of van Dantzig and van der Waerden \cite{DW},
see also pp. 46-50 of \cite{KN}, $\I(M,d_M)$ is compact in the topology
induced by $d_M$.  Moreover this theorem easy generalizes to the
following statement which we also use:

 \begin{lem}
Let $(A,d_A)$ and $(B,d_B)$ be compact metric spaces.  Let $\I(A,B)$ be
the set of isometries $f : A \to B$, and topologize it by the
compact-open topology.  Then $\I(A,B)$ is compact.
 \end{lem}

Let $\vp_k \in G$ be a sequence, and extract a subsequence converging
in $\I(M,d_M)$ to some $\vp \in \I(M,d_M)$.  Fix a Riemannian metric on
$M$ with corresponding norm $|\cdot|$.  We use the boundness of
derivatives from theorem~2 to
deduce that the sequence $(\vp_k)_\ast$ of first derivatives must be
bounded.  Otherwise we can extract a subsequence to assure that
 $$
|(\vp_k)_\ast(p_k)v_k| \to \infty
 $$
where $|v_k|=1$ and $v_k \in T_{p_k}M$.  Now extract subsequences to
assure that $p_k \to p \in M$ and $v_k \to v \in T_pM$, $|v|=1$.
Let $\xi$ be a smooth vector field extending $v$ and such that $|\xi|=1$
on a neighborhood of $p$.  By 5.4a of \cite{NW} and arguments from
chapter 1, there exists a number $r > 0$ and
a smooth family $f_q : r\D \to M$
for $q \in U$, where $U$ is a neighborhood of $p$, such that
$f_q(0) = q$ and $(f_q)_\ast(0)e = \xi_q$.  Now we obtain a sequence
of pseudoholomorphic mappings $\vp_k \circ f_{p_k} : r\D \to M$
with the property that
 $$
|(\vp_k \circ f_{p_k})_\ast(0)e| \to \infty,
 $$
which contradicts hyperbolicity.

Note also that we might obtain the boundness of the sequence of
derivatives by the nonlinear Schwarz lemma of \cite{Gr} using the
compactness.  Such a Schwarz lemma also holds for higher derivatives
\cite{Mu}, which yields the desired convergence result.

Now consider the space $J_{PH}^1(M,M)$ of 1-jets of pseudoholomorphic
mappings of $M$ to itself, i.e.\ the set of points $(p,q,\Phi)$ with
$p,q\in M$, $\Phi\in T^*_pM\otimes T_qM$, satisfying the equality
$\Phi\c J_p=J_q\c\Phi$.
This space carries a canonical almost
complex structure \cite{Gau} for which the standard projection
$\pi : J_{PH}^1(M,M) \to M \times M$ is a pseudoholomorphic mapping.
Let $\sigma : M \times M \to M$ be the projection on the second factor.
Denote the composition by $\rho = \sigma \circ \pi$.
Let $J_{PH}^1(M,M)_r$ denote the set of points $(p,q,\Phi) \in J_{PH}^1$
such that the element $\Phi \in T_p^{\ast}M \otimes T_qM$
which also could be considered as a tangent vector to $J_{PH}^1$,
satisfies $|\Phi| \leq r$.
Thus we have pseudoholomorphic mappings $\rho_r : J_{PH}^1(M,M)_r\to M$
and $\pi_r : J_{PH}^1(M,M)_r\to M\times M$. The bounded ball
$\pi_r^{-1}(p\times q)$ carries the standard complex structure and hence
is hyperbolic for every $p,\,q \in M$ and $M$ is
also hyperbolic.  Thus by theorem~3
the almost complex manifold $J_{PH}^1(M,M)_r$ is hyperbolic.

Now by the arguments above for some $r$ we have the following commutative
diagram of pseudoholomorphic mappings:

\[
\begin{array}{ccc}
& \hsb                         & \hsp   J_{PH}^{1}(M,M)_{r}   \\
& \hsb ^{j^{1}\psi}\!\nearrow  & \hsp \downarrow\! ^{\rho_{r}} \\
M \!\!\!\!\!\!\!\!\!           & \hsb \ \ \ \ \ \ \ \
\stackrel{\psi}{\longrightarrow}  & \hsp M.
\end{array}
\quad
\]

For $p,q \in M$ and an isometry $\psi \in G$ we have
 $$
d_M(p,q) \geq d_{(J_{PH}^1)_r}(j^1\psi(p),j^1\psi(q))
\geq d_M(\psi(p),\psi(q)) = d_M(p,q),
 $$
and so for every $\psi \in G$, $j^1\psi$ is an isometry with respect
to $d_M$ and $d_{(J_{PH}^1)_r}$.  Using the generalization of the
result of van Dantzig and van der Waerden quoted above, we conclude that
we may extract a subsequence to assure that $\vp_k$ converges
uniformly in $J_{PH}^1$, i.e.\ the sequence of derivatives converges
uniformly.  Exactly the same arguments above show that the derivatives
of the mappings $j^1\vp_k$ are bounded and we may consider the hyperbolic
almost complex manifold $J_{PH}^1(M,J_{PH}^1(M,M)_r)_r$ to extract a
subsequence for which we have uniform convergence in the $C^2$ sense.
Applying the argument one more time, we obtain a $\Sigma$-convergent
subsequence, thus $\vp_k \to \vp \in G$.\qed
 \end{proof}

Note that if a vector field $\xi\in\G$ is an element of the Lie
algebra of $G$ then it is complete, i.e. globally integrable.
So if also $j\xi\in\G$ then the beginning of the proof does not
use any other completeness and we obtain the following statement:

 \begin{prop}\po
No almost complex Lie group of positive dimension acts effectively
as a pseudoholomorphic transformation group on a hyperbolic
almost complex manifold.
 \end{prop}

This is an analog of the third statement of theorem~9.1 of~\cite{Ko2},
see also~\cite{Ko1}.
The first two are also valid: $\op{dim}\op{Aut}(M^{2n},J)\le 2n+n^2$
for a hyperbolic almost complex manifold with the equality iff
there is an isomorphism with the standard ball
$(M^{2n},J)\simeq (B^{2n},J_0)$ and the isotropy subgroup
$\op{Aut}(M^{2n},J)_p$ is compact for every
$p\in M$ (actually we may modify the corresponding proof that $G$ is
a closed subgroup of $\I(M)$ using the ideas from the proof of
theorem~4; we may get rid of using the lemma and
change $(J^1_{PH})_r$ to $(J^1_{PH})_C$, where $C: M\times M\to \R$ is
some bounded on compact sets function and the PH-jets $\Phi\in
T^*_pM\otimes T_qM$ satisfy $|\Phi|\le C(p,q)$).
It turns out that the statement of the proposition is still true
for another wide class of almost complex manifolds. The following
is a weak version of the general position property from~\cite{Kr}.

 \begin{dfn}
An almost complex manifold $(M^{2n},J)$ is of slightly general position
if for a dense set of points $p \in M^{2n}$ the Nijenhuis tensor
$N_J$ at every of these points satisfies
 $$
\op{Ker}[(N_J)_p(\xi, \cdot )] \neq T_{p}M^{2n}
 $$
for all vectors $\xi \in T_{p}M \backslash {\{}0{\}}$.
 \end{dfn}

 \begin{th}\po
Let an almost complex manifold $(M^{2n},J)$ be of slightly general position.
Then the Lie algebra of the automorphism group  $G = \op{Aut}(M^{2n},J)$
has the property
$\G \cap (J\G) = 0$, i.e.\ the tangent space of $G$ contains no complex
lines: $\xi \in \G \Rightarrow J\xi \not\in \G$.
 \end{th}

 \begin{proof}
If  $\xi \in \G = \op{aut}(M^{2n},J)$  then  $[\xi,J\eta] =
J[\xi,\eta]$ for any
vector field  $\eta$.  If in addition  $J\xi \in \G$  then
$[J\xi,J\eta] = J[J\xi,\eta]$.  These two equations give
$N_{J}(\xi,\eta) = 0$  at  $p$  for every  $\eta$ and hence $\xi_p=0$.
Since the set of such points $p$ is dense we have: $\xi=0$.
  \qed
 \end{proof}


\chapter{Examples}
\markboth{}{}

\hspace{13.5pt}
Numerous examples of hyperbolic and nonhyperbolic complex manifolds are
considered in the literature.  Here we consider examples with
non-integ\-ra\-ble almost complex structures.

${\bf 1^\circ}$. We start with an example of a nonhyperbolic almost
complex manifold.
The Kobayashi pseudodistance on a homogeneous almost complex manifold
is invariant (for definitions and examples of homogeneous almost complex
manifolds see \cite{KN,Y}).
Consider $S^6 = {G_2}/SU(3)$ with its well-known nonintegrable almost
complex structure $J$ that is defined by means of the octonions
(or Cayley numbers).  The definition (see \cite{KN} for the
full details) goes like this:  let
${\R}^7 = \op{Im} {\mbox{\boldmath $C\! a$}}$ be the purely
imaginary octonions and $S^6 \subset {\R}^7$ the unit
sphere.  On ${\R}^7$ there is defined a vector
product $\times$ and we define $J : T_{w}S^6 \rightarrow T_{w}S^6$ by
$\eta \mapsto {\eta}\times{w}$ where $\eta \in {\R}^7$
and $\eta \bot w$.

Let $\Ca=\H^2$ be the usual
identification of octonions as pairs of quaternions, and consider
$\H=\H\times\{0\} \subset \H^2=\Ca$ as a subspace,
$\R^3 = \op{Im}\H\subset\Ca$ as the purely imaginary quaternions and
$S^2$ as the unit sphere in $\R^3$.  Since $\R^3$ is closed
under the vector product on $\R^7$, the 2-sphere
$S^2 = S^2 \times\{0\}\subset\Ca$ is invariant under $J$ and hence
pseudoholomorphic.  Thus $S^6$ is not Brody hyperbolic and hence is not
hyperbolic.

${\bf 2^\circ}$. An almost complex structure  $J$  is called tame if
there exists a symplectic form  $\omega$  such that  $\omega(X,JX) > 0$
for any nonzero vector $X$, see \cite{Gr}.  For tame almost complex
structures on bounded domains in ${\R}^{2n}$ we have a
sufficient condition for hyperbolicity.
\pagebreak

 \begin{th}\po
Let $(D^{2n},J)$ be an almost complex domain with $J$ tame.  If there
exists a ball
$D'{\scriptstyle\supset}\!\!\!\supset\! D$ such that $J$ is the
restriction of some tame almost complex structure on $D^{\prime}$, then
$(D^{2n},J)$ is hyperbolic.
 \end{th}

 \begin{proof}
It was shown by Lafontaine \cite{Laf} that if $D$ is a
bounded domain in $({\R}^{2n},J)$, where the almost
complex structure $J$ is tamed by a symplectic form with bounded
coefficients, then $(D^{2n},J)$ is Brody hyperbolic.
We may suppose the coefficients of the taming symplectic form
on $D'$ are bounded, shrinking this ball a little otherwise.
Now consider a diffeomorphism of $D^{\prime}$ onto
${\R}^{2n}$ which maps $D$ onto a bounded domain.
The coefficients of the transported symplectic form remain
bounded. Thus we have exhibited $D$ as a bounded domain in a tame
$({\R}^{2n},J^{\prime})$, since the diffeomorphism does
not disturb the taming condition.  Now apply theorem~2
to the closure of $D$.\qed
 \end{proof}

The theorem above is the analog of the known sufficient condition for
a domain in ${\C}^n$ to be hyperbolic; that it be bounded.  For this can
be formulated as the statement that a bounded domain in
${\R}^{2n}$ with complex structure is
hyperbolic if the complex structure can be extended to the standard
complex structure on an ambient ball.  A statement of this kind is
however false for a general almost complex domain where some extra
condition such as tameness or standard integrability is not imposed.
For  $S^6 \backslash {\{}p_0{\}}$  is diffeomorphic to
${\R}^{6}$ so we can essentially impose on
${\R}^{6}$  the restriction of the almost complex
structure on  $S^6$  discussed above.  Choosing $p_0$ not to lie on the
pseudoholomorphic 2-sphere that we found, and removing a closed ball
around $p_0$, we obtain a bounded domain in  ${\R}^{6}$
(actually an open ball), such that the almost complex structure
extends to the whole of  ${\R}^{6}$. Yet this domain is not hyperbolic.

Note that the reasoning above shows that the almost complex structure
on  $S^6$  cannot be tamed even on the complement of a small closed
ball (globally it is evident since $S^6$ is not symplectic).

${\bf 3^\circ}$.
In example~2 we saw that not every bounded domain in an almost complex
$\R^{2n}$ is hyperbolic (except maybe $\R^4$). However theorem~6
implies that every point possesses a neighborhood which is hyperbolic.
Thus it follows from theorem~3 that for every almost complex manifold
$M$ with finite dimensional ($C^\infty$-)parametrized almost complex
structure $J^{(\tau)}$, $\tau\in\R^p$, the property of
$(M,J^{(\tau)})$ being hyperbolic is open with respect to (the usual)
topology on $\R^p$. This can be generalized.

The following statement gives examples of non-integrable hyperbolic
almost complex structures on compact manifolds, for example
on the product of orientable Riemann surfaces each of genus $g_i\ge2$
or on a closed Hermitian manifold with holomorphic sectional curvature
bounded above by a negative constant.

 \begin{th}\po
Let $(M,J)$ be a compact hyperbolic almost complex manifold. Then for a
small neighborhood $\ok$ of $J$ in $C^\infty$ topology every almost
complex manifold $(M,J')$ with $J'\in\ok$ is hyperbolic.
 \end{th}

 \begin{proof}
Assume that the statement is false. Then for every sequence of neighborhoods
$\ok_k$ of $J$ shrinking to $J$ there exists a sequence $J_k\in\ok_k$
of almost complex structures such that the manifold $(M,J_k)$ is not
hyperbolic. For example one may take a sequence of neighborhoods
$\ok^{(m)}_k$ to be the ball of radius $\dfrac1k$ around $J$ in
$C^m$-norm (or $W^{m,p}$-norm),
intersect it with $\ok$ and then take the diagonal subsequence. According
to theorem~2 for every $J_k$ there exists a nonconstant pseudoholomorphic
curve $f_k:\C\to(M,J_k)$. The sequence of almost complex structures $J_k$
tends to $J$ in $C^\infty$ topology. As in the proof of theorem~2 we may
apply Brody's reparametrization lemma to obtain uniform boundness of
derivatives on compact sets and the condition $|(f_k)_*(0)e|=1$,
and thus due to the compactness theorem~B.4.2.
of~\cite{MS} there exists a subsequence converging uniformly
with all derivatives on compact sets to some pseudoholomorphic
curve $f:\C\to(M,J)$. Since $|f_*(0)e|=1$ this curve is nonconstant which
contradicts hyperbolicity of $(M,J)$. \qed
 \end{proof}

%
%
%
%

%

\bigskip
\bigskip
\bigskip
\bigskip
\bigskip
\bigskip
\bigskip
\bigskip
\bigskip

\begin{itemize}

\item
{\it Address:}
{\footnotesize
 \begin{itemize}
  \item
P.O. Box 546, 119618, Moscow, Russia
  \item
Chair of Applied Mathematics, Moscow State Technological University
\linebreak
{\rm n. a.} Baumann, Moscow, Russia
 \end{itemize}
}

{\it E-mail:} \quad
{\footnotesize
lychagin\verb"@"glas.apc.org or borkru\verb"@"difgeo.math.msu.su
}

\item
{\it Address:}
{\footnotesize
 \begin{itemize}
  \item
Department of Mathematics and Statistics, University of Tromso, N-9037
Tromso, Norway
 \end{itemize}
}

{\it E-mail:} \quad
{\footnotesize
marius\verb"@"math.uit.no
}

\end{itemize}

\end{document}